\documentclass[12pt]{article}
\usepackage{amsfonts}
\usepackage{a4wide}
\usepackage{amsmath,amscd,amsthm,a4,amssymb}

\begin{document}

\begin{center}
{\Large \textbf{On the representation by bivariate ridge functions}}

\ 

\textbf{Rashid A. Aliev} \ 

Institute of Mathematics and Mechanics, NAS of Azerbaijan, Baku, Azerbaijan

Baku State University, Baku, Azerbaijan

{e-mail:} {aliyevrashid@hotmail.ru}

\ 

\textbf{Vugar E. Ismailov}\footnote{%
Corresponding author} \ 

Institute of Mathematics and Mechanics, NAS of Azerbaijan, Baku, Azerbaijan

Baku Business University, Baku, Azerbaijan

{e-mail:} {vugaris@mail.ru}
\end{center}

\textbf{Abstract.} We consider the problem of representation of a bivariate
function by sums of ridge functions. We show that if a function of a certain
smoothness class is represented by a sum of finitely many, arbitrarily
behaved ridge functions, then it can also be represented by a sum of ridge
functions of the same smoothness class. As an example, this result is
applied to a homogeneous constant coefficient partial differential equation.

\smallskip

\textbf{2010 Mathematics Subject Classification:} 26B40, 35C99, 35L25.

\smallskip

\textbf{Keywords:} Cauchy functional equation; ridge function; plane wave;
representation; smoothness.

\bigskip

\begin{center}
{\large \textbf{1. Introduction}}
\end{center}

Last 30 years have seen a growing interest in the study of special
multivariate functions called ridge functions. This interest is due to
applicability of such functions in various research areas. A \textit{ridge
function} is a multivariate function of the form%
\begin{equation*}
g(\mathbf{a\cdot x})=g(a_{1}x_{1}+\cdot \cdot \cdot +a_{n}x_{n}),
\end{equation*}%
where $g:\mathbb{R}\rightarrow \mathbb{R}$ and $\mathbf{a}=\left(
a_{1},...,a_{n}\right) $ is a fixed vector (direction) in $\mathbb{R}%
^{n}\backslash \left\{ \mathbf{0}\right\} $. These functions and their
linear combinations find applications in computerized tomography (see, e.g., 
\cite{72,97,111}), in statistics (especially, in the theory of projection
pursuit and projection regression; see, e.g., \cite{13,33}) and in the
theory of neural networks (see, e.g., \cite{Is2,M,119}). Ridge functions are
also widely used in modern approximation theory as an effective and
convenient tool for approximating complicated multivariate functions (see,
e.g., \cite{54,88,102,117}). We refer the reader to Pinkus \cite{Pinkus1}
for more on ridge functions and application areas.

It should be remarked that ridge functions have been used in the theory of
partial differential equations under the name of \textit{plane waves} (see,
e.g., \cite{69}). In general, linear combinations of ridge functions with
fixed directions occur in the study of hyperbolic constant coefficient
partial differential equations. For example, assume that $(\alpha _{i},\beta
_{i}),~i=1,...,r,$ are pairwise linearly independent vectors in $\mathbb{R}%
^{2}$. Then the general solution to the homogeneous equation

\begin{equation*}
\prod\limits_{i=1}^{r}\left( \alpha _{i}\frac{\partial }{\partial x}+\beta
_{i}\frac{\partial }{\partial y}\right) u(x,y)=0\eqno(1.1)
\end{equation*}%
are all functions of the form

\begin{equation*}
u(x,y)=\sum\limits_{i=1}^{r}v_{i}(\beta _{i}x-\alpha _{i}y)\eqno(1.2)
\end{equation*}%
for arbitrary univariate functions $v_{i}$, $i=1,...,r$, from the class $%
C^{r}(\mathbb{R}).$

Note that the solution of Eq. (1.1) is the sum of bivariate ridge functions.
Sums of bivariate ridge functions also occur in basic mathematical problems
of computerized tomography. For example, Logan and Shepp \cite{69} (the term
``ridge function" was coined by them) considered the problem of
reconstructing a given but unknown function $f(x,y)$ from its integrals
along certain lines in the plane. More precisely, let $D$ be the unit disk
in the plane and a function $f(x,y)$ be square integrable and supported on $%
D $. We are given projections $P_{f}(t,\theta )$ (integrals of $f$ along the
lines $x\cos \theta +y\sin \theta =t$) and looking for a function $g=g(x,y)$
of minimum $L_{2}$ norm, which has the same projections as $f:$ $%
P_{g}(t,\theta _{j})=P_{f}(t,\theta _{j}),$ $j=0,1,...,n-1$, where angles $%
\theta _{j}$ generate equally spaced directions, i.e. $\theta _{j}=\frac{%
j\pi }{n},$ $j=0,1,...,n-1.$ The authors of \cite{69} showed that this
problem of tomography is equivalent to the problem of $L_{2}$-approximation
of the function $f$ by sums of bivariate ridge functions with equally spaced
directions $(\cos \theta _{j},\sin \theta _{j})$, $j=0,1,...,n-1.$ They gave
a closed-form expression for the unique function $g(x,y)$ and showed that
the unique polynomial $P(x,y)$ of degree $n-1$ which best approximates $f$
in $L_{2}(D)$ is determined from the above $n$ projections of $f$ and can be
represented as a sum of $n$ bivariate ridge functions.

In this paper, we are interested in the problem of smoothness in
representation by sums of bivariate ridge functions with finitely many fixed
directions. Assume we are given $n$ pairwise linearly independent directions 
$(a_{i},b_{i}),$ $i=1,...,n,$ in $\mathbb{R}^{2}$ and a function $F:\mathbb{R%
}^{2}\rightarrow \mathbb{R}$ of the form%
\begin{equation*}
F(x,y)=\sum_{i=1}^{n}g_{i}(a_{i}x+b_{i}y).\eqno(1.3)
\end{equation*}%
Assume in addition that $F$ is of a certain smoothness class, what can we
say about the smoothness of $g_{i}$? The case $n=1$ is obvious. In this
case, if $F\in C^{k}(\mathbb{R}^{2})$, then for a vector $(c,d)\in \mathbb{R}%
^{2}$ satisfying $a_{1}c+b_{1}d=1$ we have that $g_{1}(t)=F(ct,dt\mathbf{)}$
is in $C^{k}(\mathbb{R}).$ The same argument can be carried out for the case 
$n=2.$ In this case, since the vectors $(a_{1},b_{1})$ and $(a_{2},b_{2})$
are linearly independent, there exists a vector $(c,d)\in \mathbb{R}^{2}$
satisfying $a_{1}c+b_{1}d=1$ and $a_{2}c+b_{2}d=0.$ Therefore, we obtain
that the function $g_{1}(t)=F(ct,dt\mathbf{)}-g_{2}(0)$ is in the class $%
C^{k}(\mathbb{R})$. Similarly, one can verify that $g_{2}\in C^{k}(\mathbb{R}%
)$.

The picture drastically changes if the number of directions $n\geq 3$. For $%
n=3$, there are ultimately smooth functions which decompose into sums of
very badly behaved ridge functions. This phenomena comes from the classical
Cauchy Functional Equation. This equation,%
\begin{equation*}
f(x+y)=f(x)+f(y),\text{ }f:\mathbb{R\rightarrow R}\text{,}\eqno(1.4)
\end{equation*}%
looks very simple and has a class of simple solutions $f(x)=cx,$ $c\in 
\mathbb{R}$. Nevertheless, it easily follows from the Hamel basis theory
that the Cauchy Functional Equation has also a large class of wild
solutions. These solutions are called ``wild" because they are extremely
pathological over reals. They are, for example, not continuous at a point,
not monotone at an interval, not bounded at any set of positive measure
(see, e.g., \cite{1}). Let $g$ be any wild solution of the equation (1.4).
Then the zero function can be represented as%
\begin{equation*}
0=g(x)+g(y)-g(x+y).\eqno(1.5)
\end{equation*}%
Note that the functions involved in (1.5) are bivariate ridge functions with
the directions $(1,0)$, $(0,1)$ and $(1,1)$ respectively. This example shows
that for smoothness of the representation (1.3) one must impose additional
conditions on the representing functions $g_{i},$ $i=1,...,n.$

Such additional conditions are recently found by Pinkus \cite{Pinkus}. He
proved that for a large class of representing functions $g_{i}$, the
representation is smooth. That is, if apriori assume that in the
representation (1.3), the functions $g_{i}$ belong to a certain class of
\textquotedblleft well behaved functions", then they have the same degree of
smoothness as the function $F.$ As the mentioned class of \textquotedblleft
well behaved functions" one may take, e.g., the set of functions that are
continuous at a point, bounded on one side on a set of positive measure,
monotonic at an interval, Lebesgue measurable, etc. (see \cite{Pinkus}).
Konyagin and Kuleshov \cite{Konyagin} proved that in (1.3) the functions $%
g_{i}$ inherit smoothness properties of $F$ (without additional assumptions
on $g_{i}$) if and only if the directions $\mathbf{a}^{i}$ are linearly
independent.\ Note that the results of Pinkus and also Konyagin and Kuleshov
are valid not only in bivariate but also in multivariate case.

In this paper, we study a different aspect of the problem of representation
by ridge functions. Assume in the representation (1.3) $F\in C^{k}(\mathbb{R}%
^{2})$ but the functions $g_{i}$ are arbitrary. That is, we allow very badly
behaved functions (for example, not continuous at any point). Can we write $%
F $ as a sum $\sum_{i=1}^{n}f_{i}(a_{i}x+b_{i}y)$ but with the $f_{i}\in
C^{k}(\mathbb{R)}$, $i=1,...,n$? We see that the answer to this question is
positive as expected. For the sake of convenience we state the result over $%
\mathbb{R}^{2}$, but in fact it holds over any open set in $\mathbb{R}^{2}$.

Note that the above problem is not elementary as it seems. There are cases
when representation with good functions is not possible. Such situations
happen over closed sets with no interior. In \cite{IP}, Ismailov and Pinkus
presented an example of a function of the form 
\begin{equation*}
F(x,y)=g_{1}(a_{1}x+b_{1}y)+g_{2}(a_{2}x+b_{2}y),
\end{equation*}%
that is bounded and continuous on the union of two straight lines but such
that both $g_{1}$ and $g_{2}$ are necessarily discontinuous, and thus cannot
be replaced with continuous functions $f_{1}$ and $f_{2}$.

The result of this paper can be applied to a higher order partial
differential equation in two variables if its solution is given by a sum of
sufficiently smooth plane waves (see, for example, Eq. (1.1)). Based on our
theorem below, in this case, one can demand only smoothness of the sum and
dispense with smoothness of the plane wave summands.

\bigskip

\bigskip

\begin{center}
{\large \textbf{2. Smoothness in bivariate ridge function representation}}
\end{center}

In this section we prove the following theorem.

\bigskip

\textbf{Theorem 2.1.} \textit{Assume $(a_{i},b_{i})$, $i=1,...,n$ \ are
pairwise linearly independent vectors in $\mathbb{R}^{2}$. Assume that a
function $F\in C^{k}(\mathbb{R}^{2})$ has the form}

\begin{equation*}
F(x,y)=\sum_{i=1}^{n}g_{i}(a_{i}x+b_{i}y),
\end{equation*}%
\textit{where $g_{i}$ are arbitrary univariate functions and $k\geq n-2.$
Then $F$ can be represented also in the form}

\begin{equation*}
F(x,y)=\sum_{i=1}^{n}f_{i}(a_{i}x+b_{i}y),
\end{equation*}%
\textit{where the functions $f_{i}\in C^{k}(\mathbb{R})$, $i=1,...,n$.}

\bigskip

\begin{proof} Since the vectors $(a_{n-1},b_{n-1})$ and $(a_{n},b_{n})$
are linearly independent, there is a nonsingular linear transformation $%
S:(x,y)\rightarrow (x^{^{\prime }},y^{^{\prime }})$ such that $%
S:(a_{n-1},b_{n-1})\rightarrow (1,0)$ and $S:(a_{n},b_{n})\rightarrow (0,1).$
Thus, without loss of generality we may assume that the vectors $%
(a_{n-1},b_{n-1})$ and $(a_{n},b_{n})$ coincide with the coordinate vectors $%
e_{1}=(1,0)$ and $e_{2}=(0,1)$ respectively. Therefore, to prove the theorem
it is enough to show that if a function $F\in C^{k}(\mathbb{R}^{2})$ is
expressed in the form

\begin{equation*}
F(x,y)=\sum_{i=1}^{n-2}g_{i}(a_{i}x+b_{i}y)+g_{n-1}(x)+g_{n}(y),\eqno(2.1)
\end{equation*}%
with arbitrary $g_{i}$, then there exist functions $f_{i}$ $\in C^{k}(%
\mathbb{R})$, $i=1,...,n$, such that $F$ is expressed also in the form

\begin{equation*}
F(x,y)=\sum_{i=1}^{n-2}f_{i}(a_{i}x+b_{i}y)+f_{n-1}(x)+f_{n}(y).\eqno(2.2)
\end{equation*}

By $\Delta _{l}^{(\delta )}f$ we denote the increment of a function $f$ in a
direction $l=(l_{1},l_{2}).$ That is,

\begin{equation*}
\Delta _{l}^{(\delta )}f(x,y)=f(x+l_{1}\delta ,y+l_{2}\delta )-f(x,y).
\end{equation*}%
We also use the notation $\frac{\partial f}{\partial l}$ which denotes the
derivative of $f$ in the direction $l$.

It is easy to check that the increment of a ridge function $g(ax+by)$ in a
direction perpendicular to $(a,b)$ is zero. Let $l_{1},...,l_{n-2}$ be unit
vectors perpendicular to the vectors $(a_{1},b_{1}),...,(a_{n-2},b_{n-2})$
correspondingly. Then for any set of numbers $\delta _{1},...,\delta
_{n-2}\in \mathbb{R}$ we have

\begin{equation*}
\Delta _{l_{1}}^{(\delta _{1})}\cdot \cdot \cdot \Delta _{l_{n-2}}^{(\delta
_{n-2})}F(x,y)=\Delta _{l_{1}}^{(\delta _{1})}\cdot \cdot \cdot \Delta
_{l_{n-2}}^{(\delta _{n-2})}\left[ g_{n-1}(x)+g_{n}(y)\right] .\eqno(2.3)
\end{equation*}

Denote the left hand side of (2.3) by $S(x,y).$ That is, set%
\begin{equation*}
S(x,y)\overset{def}{=}\Delta _{l_{1}}^{(\delta _{1})}\cdot \cdot \cdot
\Delta _{l_{n-2}}^{(\delta _{n-2})}F(x,y).
\end{equation*}%
Then from (2.3) it follows that for any real numbers $\delta _{n-1}$and $%
\delta _{n}$,

\begin{equation*}
\Delta _{e_{1}}^{(\delta _{n-1})}\Delta _{e_{2}}^{(\delta _{n})}S(x,y)=0,
\end{equation*}%
or in expanded form,

\begin{equation*}
S(x+\delta _{n-1},y+\delta _{n})-S(x,y+\delta _{n})-S(x+\delta
_{n-1},y)+S(0,0)=0.
\end{equation*}%
Putting in the last equality $x=y=0,$ $\delta _{n-1}=x,$ $\delta _{n}=y$, we
obtain that

\begin{equation*}
S(x,y)=S(x,0)+S(0,y)-S(0,0).
\end{equation*}%
This means that

\begin{equation*}
\Delta _{l_{1}}^{(\delta _{1})}\cdot \cdot \cdot \Delta _{l_{n-2}}^{(\delta
_{n-2})}F(x,y)=\Delta _{l_{1}}^{(\delta _{1})}\cdot \cdot \cdot \Delta
_{l_{n-2}}^{(\delta _{n-2})}F(x,0)+\Delta _{l_{1}}^{(\delta _{1})}\cdot
\cdot \cdot \Delta _{l_{n-2}}^{(\delta _{n-2})}F(0,y).
\end{equation*}

By the hypothesis of the theorem, the derivatives $\frac{\partial ^{n-2}}{%
\partial l_{1}\cdot \cdot \cdot \partial l_{n-2}}F(x,0)$ and $\frac{\partial
^{n-2}}{\partial l_{1}\cdot \cdot \cdot \partial l_{n-2}}F(0,y)$ exist.
Denote these derivatives by $h_{1,1}$ and $h_{2,1}$ respectively. Thus, it
follows from the above formula that

\begin{equation*}
\frac{\partial ^{n-2}F}{\partial l_{1}\cdot \cdot \cdot \partial l_{n-2}}%
=h_{1,1}(x)+h_{2,1}(y).\eqno(2.4)
\end{equation*}%
Note that $h_{1,1}$ and $h_{2,1}$ belong to the class $C^{k-n+2}(\mathbb{R}%
). $

By $h_{1,2}$ and $h_{2,2}$ denote the antiderivatives of $h_{1,1}$ and $%
h_{2,1}$ satisfying the condition $h_{1,2}(0)=h_{2,2}(0)=0$ and multiplied
by the numbers $1/\cos (e_{1},^{\wedge }l_{1})$ and $1/\cos (e_{2},^{\wedge
}l_{1})$ correspondingly. That is,

\begin{eqnarray*}
h_{1,2}(x) &=&\frac{1}{\cos (e_{1},^{\wedge }l_{1})}\int_{0}^{x}h_{1,1}(z)dz;
\\
h_{2,2}(y) &=&\frac{1}{\cos (e_{2},^{\wedge }l_{1})}\int_{0}^{y}h_{2,1}(z)dz.
\end{eqnarray*}%
Here $(e,^{\wedge }l)$ denotes the angle between vectors $e$ and $l$.
Obviously, the function

\begin{equation*}
F_{1}(x,y)=h_{1,2}(x)+h_{2,2}(y)
\end{equation*}%
obeys the equality

\begin{equation*}
\frac{\partial F_{1}}{\partial l_{1}}(x,y)=h_{1,1}(x)+h_{2,1}(y).\eqno(2.5)
\end{equation*}%
From (2.4) and (2.5) we obtain that

\begin{equation*}
\frac{\partial }{\partial l_{1}}\left[ \frac{\partial ^{n-3}F}{\partial
l_{2}\cdot \cdot \cdot \partial l_{n-2}}-F_{1}\right] =0.
\end{equation*}%
Hence, for some ridge function $\varphi _{1,1}(a_{1}x+b_{1}y),$

\begin{equation*}
\frac{\partial ^{n-3}F}{\partial l_{2}\cdot \cdot \cdot \partial l_{n-2}}%
(x,y)=h_{1,2}(x)+h_{2,2}(y)+\varphi _{1,1}(a_{1}x+b_{1}y).\eqno(2.6)
\end{equation*}%
Here all the functions $h_{2,1},h_{2,2}(y),\varphi _{1,1}\in C^{k-n+3}(%
\mathbb{R}).$

Set the following functions

\begin{eqnarray*}
h_{1,3}(x) &=&\frac{1}{\cos (e_{1},^{\wedge }l_{2})}\int_{0}^{x}h_{1,2}(z)dz;
\\
h_{2,3}(y) &=&\frac{1}{\cos (e_{2},^{\wedge }l_{2})}\int_{0}^{y}h_{2,2}(z)dz;
\\
\varphi _{1,2}(t) &=&\frac{1}{a_{1}\cos (e_{1},^{\wedge }l_{2})+b_{1}\cos
(e_{2},^{\wedge }l_{2})}\int_{0}^{t}\varphi _{1,1}(z)dz.
\end{eqnarray*}%
Note that the function

\begin{equation*}
F_{2}(x,y)=h_{1,3}(x)+h_{2,3}(y)+\varphi _{1,2}(a_{1}x+b_{1}y)
\end{equation*}%
obeys the equality

\begin{equation*}
\frac{\partial F_{2}}{\partial l_{2}}(x,y)=h_{1,2}(x)+h_{2,2}(y)+\varphi
_{1,1}(a_{1}x+b_{1}y).\eqno(2.7)
\end{equation*}%
From (2.6) and (2.7) it follows that

\begin{equation*}
\frac{\partial }{\partial l_{2}}\left[ \frac{\partial ^{n-4}F}{\partial
l_{3}\cdot \cdot \cdot \partial l_{n-2}}-F_{2}\right] =0.
\end{equation*}%
The last equality means that for some ridge function $\varphi
_{2,1}(a_{2}x+b_{2}y),$

\begin{equation*}
\frac{\partial ^{n-4}F}{\partial l_{3}\cdot \cdot \cdot \partial l_{n-2}}%
(x,y)=h_{1,3}(x)+h_{2,3}(y)+\varphi _{1,2}(a_{1}x+b_{1}y)+\varphi
_{2,1}(a_{2}x+b_{2}y).\eqno(2.8)
\end{equation*}%
Here all the functions $h_{1,3},$ $h_{2,3},$ $\varphi _{1,2},$ $\varphi
_{2,1}\in C^{k-n+4}(\mathbb{R}).$

Note that in the left hand sides of (2.4), (2.6) and (2.8) we have the mixed
directional derivatives of $F$ and the order of these derivatives is
decreased by one in each consecutive step. Continuing the above process,
until it reaches the function $F$, we obtain the desired result.
\end{proof}

Theorem 2.1 can be applied to Eq. (1.1) as follows.

\bigskip

\textbf{Corollary 2.2.} \textit{Assume a function $u\in C^{r}(\mathbb{R}^{2})
$ is of the form (1.2) with arbitrarily behaved $v_{i}.$ Then $u$ is a
solution to the Equation (1.1).}

\bigskip

\textbf{Remark.} Some polynomial terms appear while attempting to obtain a
smoothness result in multivariate case. In \cite{2}, we proved that if a
function $f(x_{1},...,x_{n})$ of a certain smoothness class is represented
by a sum of $r$ arbitrarily behaved ridge functions, then, under suitable
conditions, it can be represented by a sum of ridge functions of the same
smoothness class and some $n$-variable polynomial of a certain degree. The
appearance of a polynomial term is mainly related to the fact that in $%
\mathbb{R}^{n}$ ($n\geq 3)$ there are many directions orthogonal to a given
direction. Note that a polynomial term also appears in verifying if a given
function of $n$ variables ($n\geq 3$) is a sum of ridge functions (see \cite%
{25}). However, paralleling the above theorem, we conjecture that if a
multivariate function of a certain smoothness class is represented by a sum
of arbitrarily behaved ridge functions, then it can also be represented by a
sum of ridge functions of the same smoothness class.


\bigskip

\bigskip

\begin{thebibliography}{99}
\bibitem{1} J. Acz\'{e}l, \textit{Functional Equations and their Applications%
}, Academic Press, New York, 1966.

\bibitem{2} R. A. Aliev and V. E. Ismailov, On a smoothness problem in ridge
function representation, \textit{Adv. in Appl. Math.} \textbf{73} (2016),
154-169.

\bibitem{13} E. J. Cand\`{e}s, Ridgelets: estimating with ridge functions, 
\textit{Ann. Statist.} \textbf{31} (2003), 1561-1599.

\bibitem{25} P. Diaconis and M. Shahshahani, On nonlinear functions of
linear combinations, \textit{SIAM J. Sci. Stat. Comput.} \textbf{5} (1984),
175-191.

\bibitem{33} J. H. Friedman and W. Stuetzle, Projection pursuit regression, 
\textit{J. Amer. Statist. Assoc.} \textbf{76} (1981), 817-823.

\bibitem{54} V. E. Ismailov, Characterization of an extremal sum of ridge
functions, \textit{J. Comp. Appl. Math.} \textbf{205} (2007), 105-115.

\bibitem{IP} V. E. Ismailov, A. Pinkus, Interpolation on lines by ridge
functions. \textit{J. Approx. Theory} \textbf{175} (2013), 91-113.

\bibitem{Is2} V. E. Ismailov, Approximation by ridge functions and neural
networks with a bounded number of neurons, \textit{Appl. Anal.} \textbf{94}
(2015), 2245-2260.

\bibitem{69} F. John, \textit{Plane Waves and Spherical Means Applied to
Partial Differential Equations}, Interscience, New York, 1955.

\bibitem{72} I. Kazantsev, Tomographic reconstruction from arbitrary
directions using ridge functions, \textit{Inverse Problems} \textbf{14}
(1998), 635-645.

\bibitem{Konyagin} S. V. Konyagin, A. A. Kuleshov, On the continuity of
finite sums of ridge functions (Russian), \textit{Mat. Zametki} \textbf{98}
(2015), 308-309; English transl. in \textit{Math. Notes} \textbf{98} (2015),
336-338.

\bibitem{88} A. Kro\'{o}, On approximation by ridge functions, \textit{%
Constr. Approx.} \textbf{13} (1997), 447-460.

\bibitem{97} B. F. Logan and L. A. Shepp , Optimal reconstruction of a
function from its projections, \textit{Duke Math. J.} \textbf{42} (1975),
645-659.

\bibitem{102} V. E. Maiorov, On best approximation by ridge functions, 
\textit{J. Approx. Theory} \textbf{99} (1999), 68-94.

\bibitem{M} V. Maiorov and A. Pinkus, Lower bounds for approximation by MLP
neural networks, \textit{Neurocomputing} \textbf{25} (1999), 81--91.

\bibitem{111} F. Natterer, \textit{The Mathematics of Computerized Tomography%
}, Wiley, New York, 1986.

\bibitem{117} P. P. Petrushev, Approximation by ridge functions and neural
networks, \textit{SIAM J. Math. Anal.} \textbf{30} (1998), 155-189.

\bibitem{Pinkus1} A. Pinkus, Ridge Functions, Cambridge Tracts in
Mathematics, 205. Cambridge University Press, Cambridge, 2015.

\bibitem{Pinkus} A. Pinkus, Smoothness and uniqueness in ridge function
representation, \textit{Indag. Math.} (N.S.) \textbf{24} (2013), no. 4,
725--738.

\bibitem{119} A. Pinkus, Approximation theory of the MLP model in neural
networks, \textit{Acta Numerica} \textbf{8} (1999), 143-195.
\end{thebibliography}
\end{document}